\documentclass{article}

\usepackage{arxiv}

\usepackage[utf8]{inputenc} 
\usepackage[T1]{fontenc}    
\usepackage{hyperref}       
\usepackage{url}            
\usepackage{booktabs}       
\usepackage{amsmath,amssymb,amsfonts}     
\usepackage{nicefrac}       
\usepackage{microtype}      
\usepackage{lipsum}
\usepackage{graphics}
\usepackage{graphicx}
\usepackage{mathrsfs}
\usepackage{multirow}
\usepackage{enumerate}
\usepackage{hyperref}
\usepackage{cleveref}

\title{Parallel Randomized Algorithm for Chance Constrained Program}

\author{
  Xun Shen \\
  Department of Statistical Sciences\\
  The Graduate University for Advanced Studies\\
  Tokyo, Japan 106-8569\\
  \texttt{shenxun@ism.ac.jp} \\
   \And
 Jiancang Zhuang \\
  Institute of Statistical Mathematics\\
  Research Organization of Information and Systems\\
  Tokyo, Japan 190-8562\\
  \texttt{zhuangjc@ism.ac.jp} \\
  \And
  Xingguo Zhang \\
  Department of Mechanical Systems Engineering\\
  Tokyo University of Agriculture and Technology\\
  Tokyo, Japan, 184-8588\\
  \texttt{xgzhang@go.tuat.ac.jp} \\
}

\begin{document}
\maketitle

\begin{abstract}
Chance constrained program is computationally intractable due to the existence of chance constraints, which are randomly disturbed and should be satisfied with a probability. This paper proposes a two-layer randomized algorithm to address chance constrained program. Randomized optimization is applied to search the optimizer which satisfies chance constraints in a framework of parallel algorithm. Firstly, multiple decision samples are extracted uniformly in the decision domain without considering the chance constraints. Then, in the second sampling layer, violation probabilities of all the extracted decision samples are checked by extracting the disturbance samples and calculating the corresponding violation probabilities. The decision samples with violation probabilities higher than the required level are discarded. The minimizer of the cost function among the remained feasible decision samples are used to update optimizer iteratively.  Numerical simulations are implemented to validate the proposed method for non-convex problems comparing with scenario approach. The proposed method exhibits better robustness in finding probabilistic feasible optimizer.
\end{abstract}


\section{Introduction}
\label{sec:introduction}
Chance constrained program, also identified as probabilistic constrained program\cite{Prekopa}, was firstly considered in the fields of management and economics\cite{Charnes}. Chance constrained program involves constraints disturbed by random variables. These randomly disturbed constraints are required to be satisfied with specified probability levels. Thus, they are called chance constraints. In recent years, chance constrained program has been applied to automotive powertrain control, model predictive control, system identification, machine learning and many other fields\cite{Shen,CampiBook}. 

However, there is challenge to apply chance constrained program. The existence of chance constraints makes it a NP hard problem. This motivates the use of approximate approach to target chance constrained problem, and the main stream have converged to scenario or sample-based approaches \cite{Calafiore,Luedtke} in which chance constraints are replaced by deterministic constraints imposed for finite sets of independent samples extracted of the uncertain parameters. If the assumption that the constraints are convex in the decision variables for all extracted samples holds, scenario approach preserves the convexity of the reformed deterministic problem. Furthermore, the solution of the deterministic problem satisfies chance constraints in original problem with a determined bounds of probability. The bounds of probability is related to the number of the extracted samples \cite{Calafiore}. Afterwards, since the tight bounds of samples number saves the computation time, scenario approach with tight confidence bounds has been popular. In scenario approach with tight confidence bounds, the core is to improve confidence of probability on satisfying the chance constraints for the specific sample number. For instance, \cite{CampiSampling} proposed a discarding method, in which a certain proportion of parameter samples to define a set of sampled constraints and discard the rest. Therefore, chance constraints can be approximated with tight violation probabilities for fixed sample number. Moreover, a tighter bounds is provided by a repetitive scenario approach in which a priori probabilities of violation probability is utilized\cite{Cannon}. However, scenario approach still has a fatal drawback. It cannot ensure that the obtained solution is the optimal one in the probabilistic feasible domain. when the sample number becomes larger, the obtained solution becomes more conservative and finally converges to the totally robust solution which is feasible for all uncertainty realizations. Moreover, the optimal solution is highly depended on the extracted samples, thus, it is not able to ensure that the solution converges to the optimal one which satisfies the chance constraints.  

Bayesian optimization framework has been applied to optimization under unknown constraints recently\cite{Gramacy,Picheny}, which is essentially a data-driven approach for approximating the optimizer of the program. Stimulated by \cite{Gramacy}, this paper addresses chance constrained program with a two-layer sample-based numerical approach. Basic randomized optimization is refined to a parallel sampling algorithm to search the optimizer in the probabilistic feasible domain in which chance constraints is satisfied. In the first layer, Multiple decision samples are extracted uniformly in the decision domain beyond the consideration of chance constraints. Then, verifications of the violation probability for extracted decision samples are implemented and the ones with violation probabilities higher than the required level are discarded. The minimizer of the cost function among the remained feasible decision samples are used to update optimizer iteratively. 

The rest of the paper is organised as follows. Section \ref{sec:background} gives brief background of chance constrained program formerly and then present the problem description. In Section \ref{sec:proposed}, proposed algorithms is introduced after brief discussion of the randomized optimization. The numerical simulation for validating the sampling algorithm is presented in Section \ref{sec:numerical}, using a non-convex program with chance constraints as targeted problem. Finally, Section \ref{sec:conclusion} concludes the whole study. 

\section{Background and Problem Description}
\label{sec:background}
\begin{figure}
\centering
\includegraphics[scale=0.6]{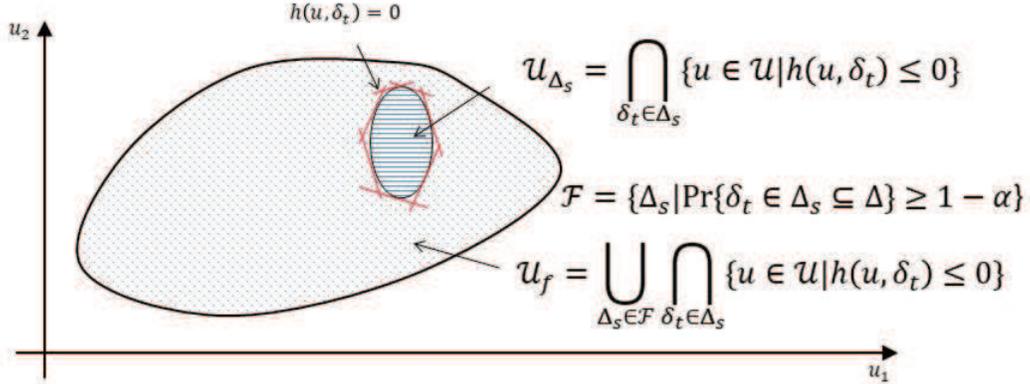}
\centering
\caption{Formulation of probabilistic feasible domain.}
\label{prob_fea_dom}
\end{figure}

Chance constrained program can be generally expressed as:
\begin{equation}
\label{eq_ccp}
\begin{split}
&\underset{u\in\mathscr{U}}{\text{min}} \,\, J(u) \\
&s.t.\quad  \text{Pr}\{h(u,\delta)\leq 0\}\geq 1-\alpha,\ \delta\in\Delta,\ \alpha\in (0,1) 
\end{split}
\end{equation} 
where $u\in\mathscr{U}\subset\mathbb{R}^{n_u}$ is the decision variable, the decision variable domain $\mathscr{U}$ is bounded, $\delta\in\Delta\subset\mathbb{R}^{n_{\delta}}$ is an uncertain parameter, the set $\Delta$ is the sample space of $\delta$ and $\text{Pr}$ is a probability measure defined on $\Delta$, $\alpha$ is a given probability level. Moreover, $J(u):\mathscr{U}\times\Delta\ \rightarrow\ \mathbb{R}$ and $\forall\delta\in\Delta,h(u,\delta):\mathscr{U}\times\Delta\ \rightarrow\ \mathbb{R}$ are continuous and differentiable in $u$. Chance constraints appearing in Eq. \eqref{eq_ccp} emerge in various applications and can be regarded as a compromise of hard constraints. Hard constraints require satisfying the constraints $h(u,\delta)\leq 0$ for all values $\delta\in\Delta$ which could be too costly and even impossible. 

The feasible decision domain can be denoted as $\mathscr{U}_f$. However, it is difficult to obtain the exact expression of $\mathscr{U}_f$. Denoting $\Delta_{s}\subset\Delta$ and the probability measure of $\Delta_{s}$ satisfies 
\begin{equation}
\label{eq_pmds}
\text{Pr}\{\Delta_{s}\}\geq 1-\alpha.
\end{equation} 
The feasible domain for $\delta_t\in\Delta_{s}$ is defined as
\begin{equation}
\label{eq_fds}
\mathscr{U}_{\delta_t}=\{u\in\mathscr{U}|h(u,\delta_t)\leq 0\}.
\end{equation}
Then, the feasible domain for $\Delta_{s}$ is intersection of $\mathscr{U}_{\delta_t}$ for all $\delta_t\in\Delta_{s}$, which is written as  
\begin{equation}
\label{eq_fds}
\mathscr{U}_{\Delta_s}=\bigcap_{\delta_t\in\Delta_{s}}\mathscr{U}_{\delta_t}.
\end{equation}
Considering a family of $\Delta_s$ denoted as
\begin{equation}
\label{eq_family}
\mathscr{F}=\{\Delta_{s}\subset\Delta|\text{Pr}\{\Delta_{s}\}\geq 1-\alpha\},
\end{equation}
the feasible decision domain for Eq. \eqref{eq_ccp} can be defined as:
\begin{equation}
\label{eq_family}
\mathscr{U}_{f}=\bigcup_{\Delta_{s}\in\mathscr{F}}\mathscr{U}_{\Delta_s}=\bigcup_{\Delta_{s}\in\mathscr{F}}\bigcap_{\delta_t\in\Delta_{s}}\mathscr{U}_{\delta_t}.
\end{equation}
The simple concept of the above process is illustrated in Fig. \ref{prob_fea_dom}. Obviously, even if $\mathscr{U}_{\delta_t}$ is known, it is impossible to obtain explicit $\mathscr{U}_{f}$ due to infinite times' operation of intersection and union. Thus, program expressed in \eqref{eq_ccp} is NP hard due to the chance constraint. To address program in \eqref{eq_ccp}, the following issues should be considered:
\begin{itemize}
\item How to approximate chance constraints, namely approximate the probabilistic feasible domain of decision variable;
\item How to approximate the optimizer in the probabilistic feasible domain. 
\end{itemize}

\section{Proposed Method}
\label{sec:proposed}
This section discusses the parallel optimizer exploration approach. Firstly, the structure of the proposed method is presented. Then, randomized optimization is briefly introduced. Finally, the exploration algorithm is introduced in details.

\subsection{Overview}
\begin{figure}
\centering
\includegraphics[scale=0.6]{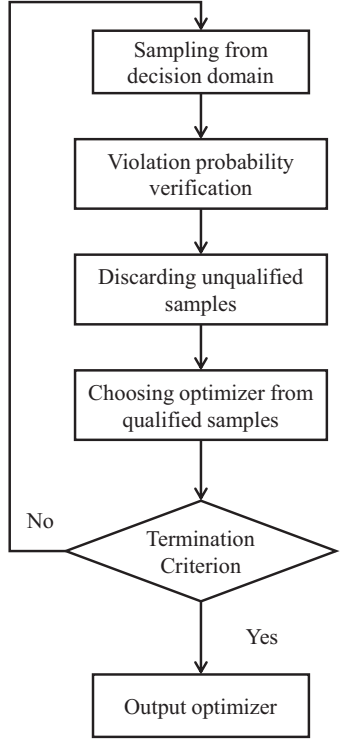}
\centering
\caption{Flow diagram of randomized optimization for chance constrained program.}
\label{ran_fra}
\end{figure}

Randomized optimization can be used to solve program in \eqref{eq_ccp}. The flow diagram of randomized optimization applied for chance constrained program is illustrated in Fig. \ref{ran_fra}. Firslty, decision samples are randomly extracted from $\mathscr{U}$ which might not satisfy the chance constraints. Thus, in the next step, violation probability of extracted samples are verified and the unqualified samples are discarded. The optimizer are chosen from the qualified samples. The process is repeated until termination criterion is satisfied and optimizer is outputted finally.

\subsection{Randomized Optimization Algorithm}
Randomized optimization algorithm is numerical optimization algorithm which does not require the gradient of the problem and it can therefore be used on functions that are not continuous or differentiable\cite{Matyas}. Moreover, if the problem is constrained and the constraints might be non-convex, randomized optimization algorithm still ensures the convergence to the optimizer \cite{Baba}. Random optimization moves to better positions iteratively in the exploration domain (the exploration domain should be a subset of feasible domain), sampling near the current position according to normal distribution or uniform distribution\cite{Dorea}.

Let $F:\mathscr{V}\subset\mathbb{R}^n\rightarrow\mathbb{R}$ be the cost function to be minimized in the feasible domain $\mathscr{V}$. Let $v\in\mathscr{V}$ denote a position or candidate optimizer in the feasible domain. The typical randomized optimization algorithm can then be denoted as algorithm 1 and expressed as following\cite{Solis}:
\begin{itemize} 
\item Step 1: Initialize $v\in\mathscr{V}$ randomly  \\
\item Step 2: Extract new sample $v'\in\mathscr{V}_v$ obeying uniform distribution or normal distribution. Here, $\mathscr{V}_v$ is a neighbour of $v$, for instance, as $
\mathscr{V}_v=\{v_f\in\mathscr{V}|\left\|v_f-v\right\|^2<\epsilon_v\}$ \\
\item Step 3: If $F(v')<F(v)$, set $v=v'$ \\
\item Step 4: Examine whether the termination criterion is met(e.g. number of iterations), if termination criterion is met, go to Step 5, otherwise, go to Step 2 and repeat \\
\item Step 5: Set $v$ as optimizer
\end{itemize}

\subsection{Optimizer Exploration Algorithm}
\begin{figure*}
\centering
\includegraphics[scale=0.4]{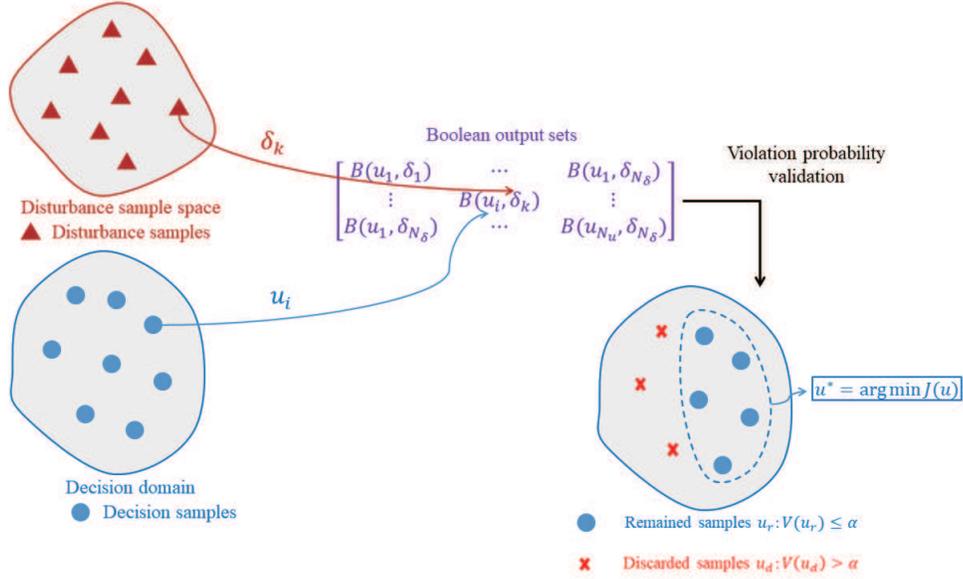}
\centering
\caption{Brief structure of optimizer exploring algorithm.}
\label{exp_opt}
\end{figure*}

The proposed optimizer exploration algorithm firstly extracts $N_u$ samples, $\{u_{1},...,u_{N_u}\}\in\mathscr{U}$ , from the decision domain. However, not all the samples satisfy the chance constraints. The violation probabilities of decision samples should be validated. Violation probability of decision $u$ is denoted as 
\begin{equation}
\label{eq_vioprob}
V(u)=\text{Pr}\{h(u,\delta)>0\},\delta\in\Delta.
\end{equation}
For given decision $u_i$ and disturbance $\delta_k$, a boolean function can be denoted as 
\begin{equation}
\label{eq_bool}
B(u_i,\delta_k)=
\left\{
\begin{split} 
&1, \text{if}\ h(u_i,\delta_k)>0 \\
&0, \text{if}\ h(u_i,\delta_k)\leq0 
\end{split}
\right.
\end{equation}
If $N_{\delta}$ samples, $\{\delta_{1},...,\delta_{N_{\delta}}\}$, are extracted from $\Delta$ independently obeying the identical distribution, the estimate violation probability can be written as 
\begin{equation}
\label{eq_vpest}
\hat{V}(u_i) = \frac{\sum_{k=1}^{N_{\delta}}B(u_i,\delta_k)}{N_{\delta}}.
\end{equation}
According to law of larger numbers, as $N_{\delta}\rightarrow\infty$, $\hat{V}(u_i)$ converges to $V(u_i)$ with probability one\cite{Feller}. Furthermore, for fixed sample number $N_{\delta}$, the mean value of $\hat{V}(u_i)$, denoted as $\mathbb{E}\{\hat{V}(u_i)\}$, equals to $V(u_i)$ and the variance $\sigma\{\hat{V}(u_i)\}$ equals to $\frac{\sigma\{V(u_i)\}}{N_{\delta}^2}$. Therefore, based on the above disturbance sampling, $V(u_{i})=Pr\{h(u_{i},\delta)>0\}, \forall i\in\{1,...,N_u\}$ can be calculated and the ones with violation probability larger than $\alpha-\alpha_{\epsilon}$ are discarded, where $\alpha_{\epsilon}$ is within $(0,\alpha)$. The boundary for discarding samples is chosen slightly smaller than the required one for a techinical consideration. Since the estimate $\hat{V}(u_i)$ is centered at $V(u_i)$ with variance, half of the estimates might be located in the position where the real violation probability is larger than the required one. By chosen a slightly smaller bound, the percentage of points located out of feasible area decreases as $N_{\delta}$ increases. The remained feasible samples are used to calculate the corresponding cost function values. The one with minimal cost value is chosen and compared with the current optimizer. If it is better than the current optimizer, it will be regarded as optimizer instead. The above steps are repeated until termination criterion is met, for instance, number of iterations. The algorithm for exploring the feasible optimizer is summarized in Algorithm 2 
\begin{itemize}
\item Step 1: Initialize optimizer $u_f^{*}\in\mathscr{U}$ randomly\\
\item Step 2: Extract $\{u_{1},,...,u_{i},...,u_{N_u}\}\in\mathscr{U}$ \\
\item Step 3: Extract $\{\delta_{1},...,\delta_{N_{\delta}}\}\in\Delta$ randomly\\
\item Step 4: Calculate approximately $V(u_{i}), \forall i\in\{1,...,N_u\}$ according to Eq. \eqref{eq_bool} and \eqref{eq_vpest}, discard all $u_{d}$ that $V(u_{d})>\alpha-\alpha_{\epsilon}$, and the remained solution set is $\mathbb{U}_f$ \\
\item Step 5: If $\mathbb{U}_f=\emptyset$, return to Step 1, else if $\mathbb{U}_f\neq\emptyset$, go to Step 6\\
\item Step 6: Set optimizer as $u^{*}=\underset{u\in\mathbb{U}_f}{\text{argmin}} J(u)$\\
\item Step 7: Replace optimizer $u_f^{*}=u^{*}$ and go to Step 8 if $J(u^{*})<J(u_f^{*})$; Return to Step 2 directly if $J(u^{*})\geq J(u_f^{*})$  \\
\item Step 8: Examine whether the termination criterion is met (e.g. number of iterations), if termination criterion is met, go to Step 6, otherwise, go to Step 1 and repeat \\
\item Step 9: Set $u_f^{*}$ as optimizer.
\end{itemize} 
Moreover, the above algorithm is briefly illustrated in Fig. \ref{exp_opt} as well. Apparently, the algorithm is parallel compared to the original randomized algorithm. It is able to applied in the hardware which supports parallel algorithm and can improve computation efficiency. 

\section{Numerical Example}
\label{sec:numerical}

\begin{figure}
\centering
\includegraphics[scale=0.6]{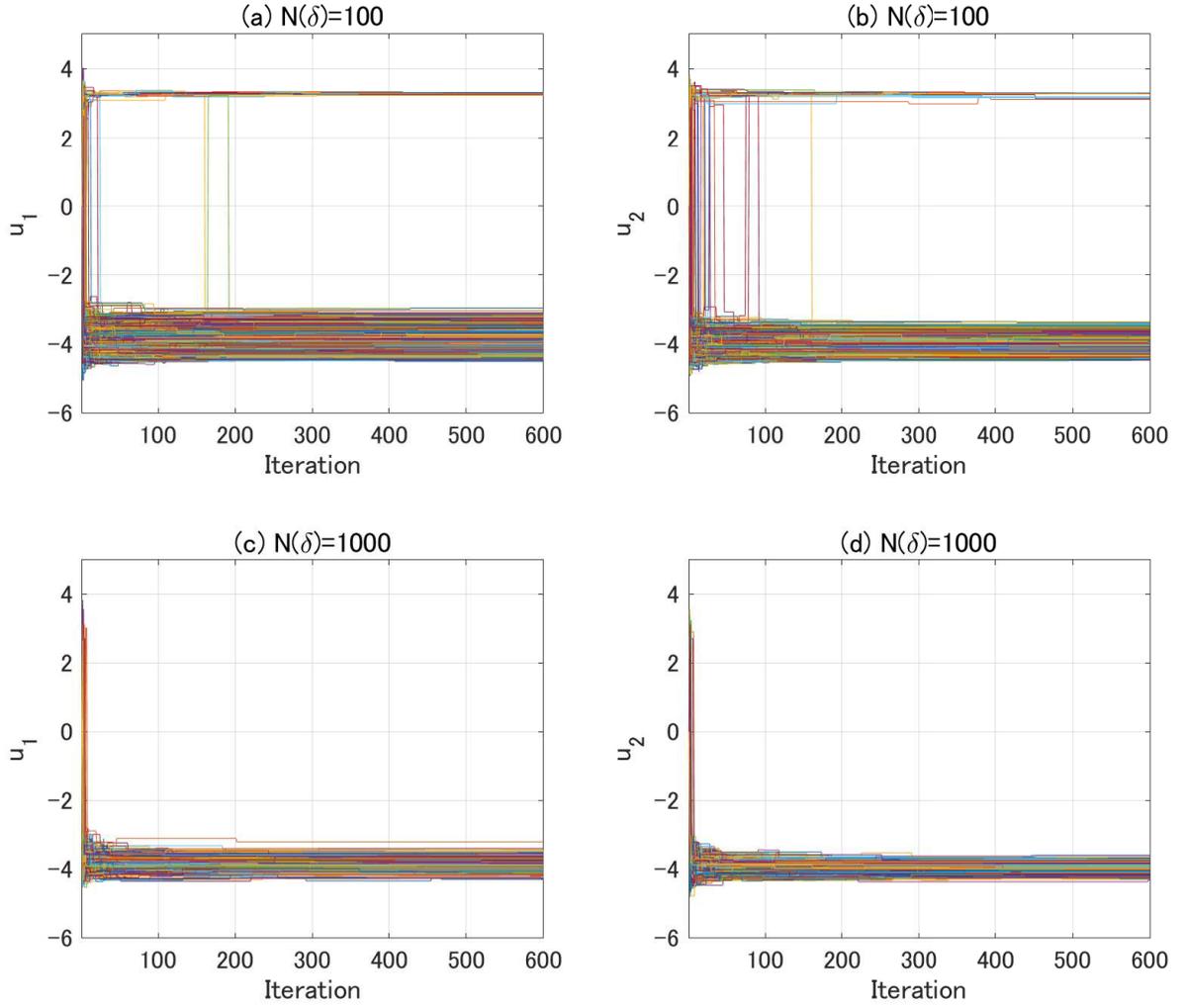}
\centering
\caption{Evolutions of the solution during computation process for proposed method (simulation of 500 times).}
\label{u_evo}
\end{figure}

\begin{figure}
\centering
\includegraphics[scale=0.6]{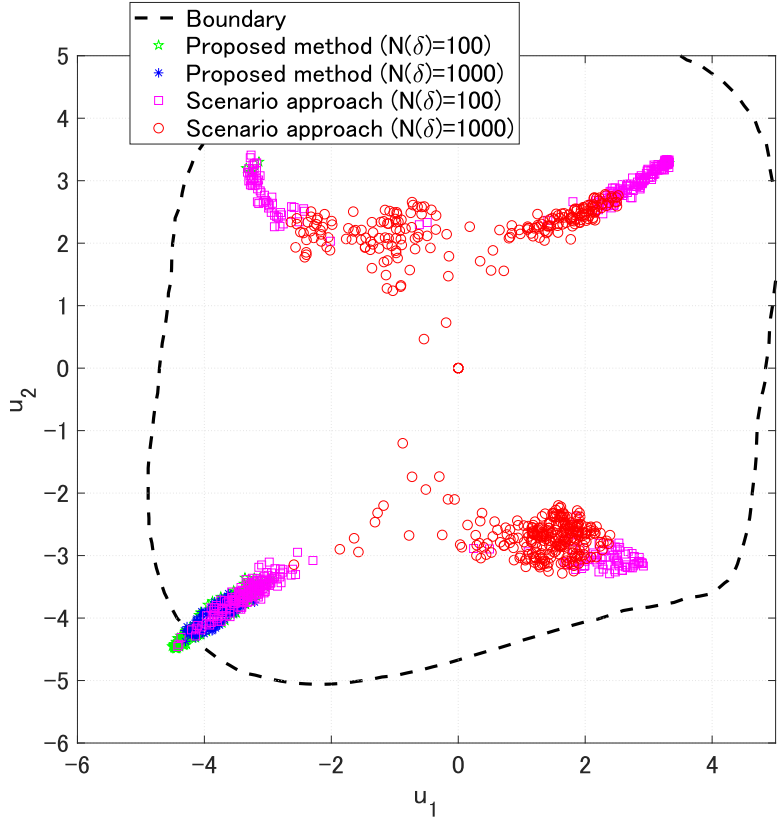}
\centering
\caption{Solution distribution (simulation of 500 times.}
\label{Sol_dis}
\end{figure}

\begin{figure}
\centering
\includegraphics[scale=0.6]{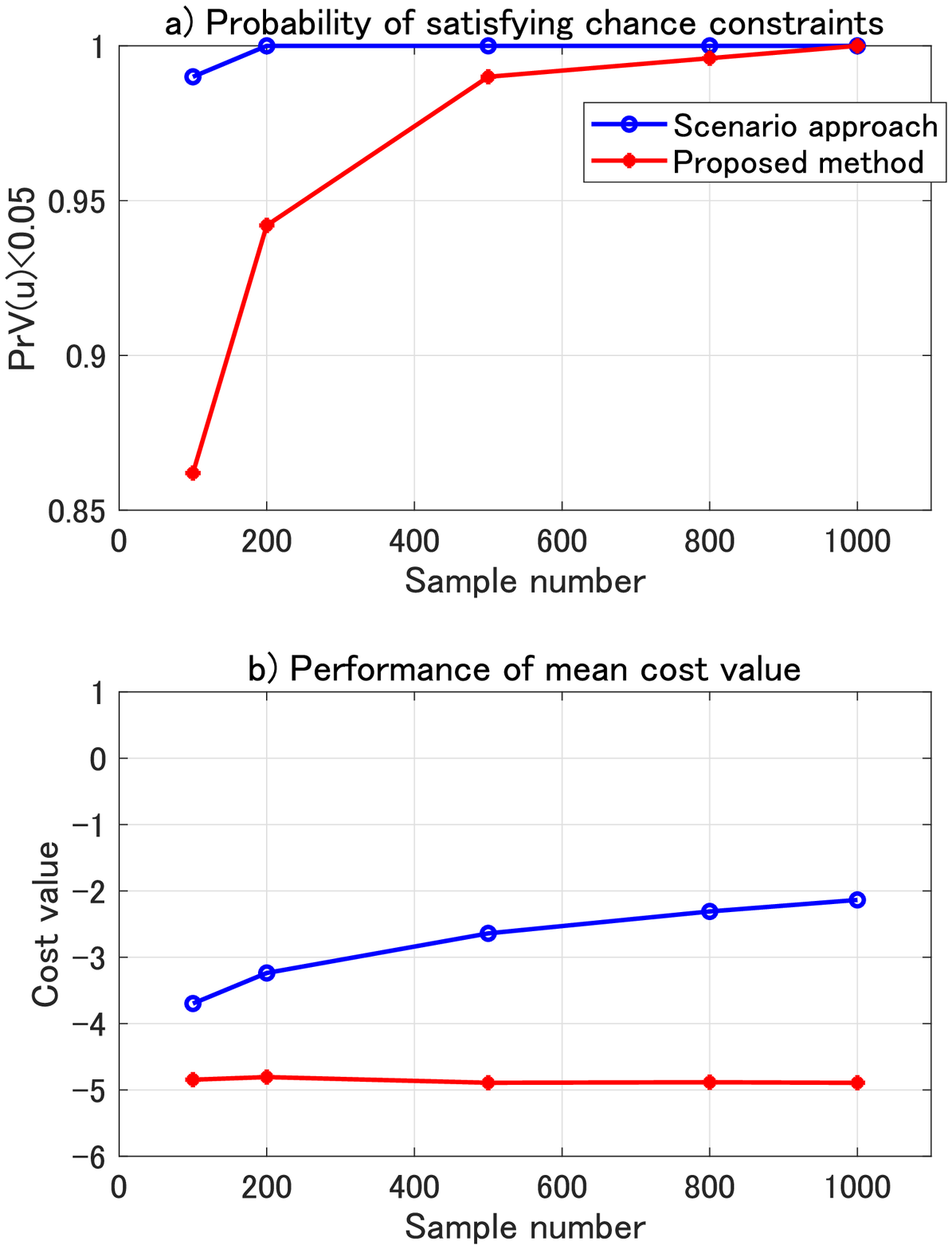}
\centering
\caption{Results of the cost values and violation probabilities (simulation of 500 times).}
\label{res_plot}
\end{figure}

This section presents numerical simulation for verifying the proposed method, comparing with scenario approach. The targeted problem is briefly introduced firstly. Also, the basic concept of scenario approach is introduced shortly which is enough for helping understand the comparison result. Furthermore, the results of optimizer exploration are presented. 

\subsection{Targeted Problem}

The targeted problem in the numerical simulation is a non-convex program with chance constraints. The decision domain is $\mathscr{U}=[-6,5]^2$. The cost function is 
\begin{equation}
\label{eq_tp_cost}
J(u) = \frac{\sum_{i=1}^2(u_i+0.5)^4-30u_i^2-20u_i}{100}.
\end{equation}
The constrained function is
\begin{equation}
\label{eq_tp_cons}
h(u,\delta) = (\sum_{i=1}^2 0.05*(u_i-a_{i}\delta)^4-b_{i}*(u_i-a_{i}\delta)^2)-(1-0.1\delta)^2,
\end{equation}
where $a_1=1.5,a_2=2,b_1=2, b_2=3$, $\delta$ is random variable which obeys normal distribution $\mathscr{N}(0,1)$. Moreover, the violation probability level is $\alpha=0.05$.

The proposed method is applied to solve the targeted problem, comparing with scenario approach presented in \cite{CampiSampling}. In scenario approach, independent samples $\delta^{(i)},i=1,...,N$ is identically extracted from $\Delta$ randomly, a deterministic convex optimization problem can be formed as
\begin{equation}
\label{eq_dcop}
\begin{split}
&\underset{u\in\mathscr{U}}{\text{min}} \,\, J(u) \\
&s.t.\quad  h(u,\delta^{(i)})\leq 0,\ i=1,...,N
\end{split}
\end{equation}
which is a standard finitely constrained optimization problem. The optimal solution $\hat{u}_N$ of the program Eq. \eqref{eq_dcop} is called as the scenario solution for program Eq. \eqref{eq_ccp} generally. Moreover, since the extractions $\delta^{(i)},i=1,...,N$ is randomly chosen, the optimal solution $\hat{u}_N$ is random variable. If $\hat{u}_N$ is expected to satisfy
\begin{equation}
\label{eq_sa_pr}
\text{Pr}^{N}\{(\delta^{(1)},...,\delta^{(N)}\in\Delta^N:V(\hat{u}_N)\leq\alpha\}\geq1-\beta, \beta\in(0,1),
\end{equation}
then, $N$ should have a lower limitation $N_l$
\begin{equation}
\label{eq_sa_n}
N\geq \frac{2}{\alpha}\text{ln}\frac{1}{\beta}+2n_u+\frac{2n_u}{\alpha}\text{ln}\frac{2}{\alpha}.
\end{equation}
Note that $\beta$ is an important factor and choosing $\beta=0$ makes $N_l=\infty$. Namely, if the number of chosen samples gets larger, the probability of satisfying the original chance constraints approaches 1. Actually, when number of chosen samples becomes infinity, the samples cover the whole sample space, not only the feasible area determined by chance constraints, it becomes totally robust.

\subsection{Optimizer Exploration}

The results of optimizer exploration are presented here. In this numerical simulation, both proposed method and scenario approach are used to solve the targeted problem. The numerical simulation adopts Monte Carlo methodology, namely, the procedure of exploring the optimizers is repeated for a large number of times and the samples of $u$ and $\delta$ are extracted obeying the identical distribution in every simulation time. 500 times of simulation were done in this validation. Moreover, for the proposed method, the parameter $\alpha_{\epsilon}$ is chosen as 0.005.

In Fig. \ref{u_evo}, evolutions of the solution during computation process for proposed method are plotted, for both cases $N(\delta)=100$ and $N(\delta)=1000$. In every simulation, the solutions converge to the small intervals for both $u_1$ and $u_2$. The final converged values of decision are plotted in Fig. \ref{Sol_dis}. Green star marks stand for solutions from proposed method using 100 disturbance samples, blue points are solutions from proposed method using 1000 disturbance samples, magenta squares are solutions from scenario approach using 100 disturbance samples and red circles present the solutions from scenario approach using 1000 disturbance samples. Apparently, proposed method can converges to the $5\%$ boundary where the cost is minimal in the probabilistic domain. However, for scenario approach, the solution cannot converge to the boundary. In every simulation, the sampled $\delta$ are different, the constraints of problem expressed by Eq. \eqref{eq_dcop} are different, the solution are therefore different and distributed in a larger area. When $N(\delta)$ is small, the considered deterministic constraints are not enough to ensure the violation probability, the solutions are mostly distributed outside the feasible area as shown by magenta squares in Fig. \ref{Sol_dis}. While, as $N(\delta)$ gets larger, the constraints become more and the solutions are consequently conservative and distributed inside the boundary as red circles in Fig. \ref{Sol_dis}. 

Furthermore, Fig. \ref{res_plot} shows the statistical results of cost values and violation probabilities in all simulations. From these results, obviously, the proposed method achieves the trade-off between violation probability and optimization in a more stable way. While, scenario approach fails to find the balance. It gets minimal cost with loss of violation probability performance in less disturbance samples' simulations. On the other hand, it provides too conservative results on violation probability with worse cost performance. The mean costs of the proposed method are all near the boundary even the samples of disturbance become more. The violation of chance constraints becomes less when there are larger sample numbers of disturbance.

\section{Conclusion and Future work}
\label{sec:conclusion}

This paper has introduced a statistical approach to chance constrained optimization. The novel idea is to approximate probabilistic feasible optimizer with a parallel randomized algorithm. The proposed method is validated by numerical simulation compared with scenario approach and exhibits better robustness on both exploring optimizer and satisfying violation probability. while, there still remains future works to be done for improving the proposed method. The current randomized optimization algorithm is totally a random one which can be improved to converge to the optimizer in less iterations. Moreover, the accuracy of probability estimate is related to the sample number of random disturbance. The quantitative analysis should be implemented to investigate the relationship between accuracy and sample number.

\bibliographystyle{unsrt}  


\end{document}